\documentclass[12pt]{article}
\usepackage{amssymb,amsmath}

\begin{document}

\baselineskip15pt

\font\fett=cmmib10
\font\fetts=cmmib7
\font\bigbf=cmbx10 scaled \magstep2
\font\bigrm=cmr10 scaled \magstep2
\font\bbigrm=cmr10 scaled \magstep3
\def\cnl{\centerline}
\def\etal{{\it et al.}}
\def\text{{}}
\def\un{^{(n)}}
\def\ut{^{(2)}}
\def\ff{{\cal F}}
\def\Ref#1{(\ref{#1})}
\def\tw{{\tilde w}}

\def\Blm{\left|}
\def\Brm{\right|}
\def\Bl{\left(}
\def\Br{\right)}
\def\nti{n\to\infty}
\def\lnti{\lim_{\nti}}
\def\law{{\cal L}}
\def\sji{\sum_{j\ge1}}
\def\Cal{\cal}

\def\real{\text{\rm I\kern-2pt R}}
\def\re{\real}
\def\expec{\text{\rm I\kern-2pt E}}
\def\ex{\expec}
\def\prob{\text{\rm I\kern-2pt P}}
\def\pr{\prob}
\def\rat{\text{\rm Q\kern-5.5pt\vrule height7pt depth-1pt\kern4.5pt}}
\def\comp{\text{\rm C\kern-4.7pt\vrule height7pt depth-1pt\kern4.5pt}}
\def\nat{\text{\rm I\kern-2pt N}}
\def\integ{{\bf Z}}
\def\qedbox{\vcenter{\hrule height.5pt\hbox{\vrule width.5pt height8pt
\kern8pt\vrule width.5pt}\hrule height.5pt}}
\def\half{{\textstyle {1 \over 2}}}
\def\quarter{{\textstyle {1 \over 4}}}
%\mathchardef\scrl="024C
%\mathchardef\scrd="0244
\def\scrl{{\Cal L}}
\def\scrd{{\Cal D}}
\def\tod{\buildrel \scrd \over \longrightarrow}
\def\eqd{\buildrel \scrd \over =}
\def\tolone{\buildrel L_1 \over \longrightarrow}
\def\toltwo{\buildrel L_2 \over \longrightarrow}
\def\tolp{\buildrel L_p \over \longrightarrow}
\def\l{\lambda}
\def\L{\Lambda}
\def\a{\alpha}
\def\b{\beta}
\def\g{\gamma}
\def\G{\Gamma}
\def\d{\delta}
\def\D{\Delta}
\def\e{\varepsilon}
\def\h{\eta}
\def\z{\zeta}
\def\th{\theta}
\def\k{\kappa}
\def\m{\mu}
\def\n{\nu}
\def\p{\pi}
\def\r{\rho}
\def\s{\sigma}
\def\S{\Sigma}
\def\t{\tau}
\def\f{\varphi}
\def\ch{\chi}
\def\ps{\psi}
\def\o{\omega}
\def\lee{\,\le\,}
\def\leee{\quad\le\quad}
\def\gee{\,\ge\,}
\def\geee{\quad\ge\quad}
\def\scrn{{\Cal N}}
\def\scra{{\Cal A}}
\def\scrf{{\Cal F}}
\def\var{\text{\rm Var\,}}
\def\cov{\text{\rm Cov\,}}

\def\sin{\sum_{i=1}^n}
\def\sn{\sum_{i=1}^n}
\def\cross{\times}

\def\nin{\noindent}
\def\bsk{\bigskip}
\def\msk{\medskip}
\def\widebar{\bar}

\def\sno{\sum_{n\ge 0}}
\def\sj{\sum_{j\ge 0}}
\def\proof{\noindent{\bf Proof.}\ }
\def\remark #1. {\medbreak\noindent{\bf Remark #1.}\quad}
\def\remarks{\medbreak\noindent{\bf Remarks.}\quad}
\def\head #1/ {\noindent{\bf #1.}\medskip}
\def\dtv{d_{TV}}
\def\Po{\text{\rm Po\,}}
\def\Bi{\text{\rm Bi\,}}
\def\Be{\text{\rm Be\,}}
\def\CP{\text{\rm CP\,}}

\def\pb{\hbox{{\fett\char'031}}}
\def\pbh{\text{\bf{\hbox{\^{\fett\char'031}}}}}
\def\bp{\pb}
\def\lb{\hbox{{\fett\char'025}}}
\def\lbs{\hbox{{\fetts\char'025}}}
\def\bl{\lb}
\def\rb{\hbox{\fett\char'032}}
\def\sb{\hbox{\fett\char'033}}
\def\iid{\text{\rm independent\ and\ identically\ distributed}}

\def\and{\ \text{\rm{ and }}\ }
\def\for{\ \text{\rm{ for }}\ }
\def\forall{\ \text{\rm{ for all }}\ }
\def\forsome{\ \text{\rm{ for some }}\ }
\def\tif{\text{{\rm if\ }}}
\def\tin{\text{{\rm in\ }}}
\def\lcm{\text{\rm{l.c.m.\,}}}

\def\xx{{\cal X}}
\def\supp{{\rm supp\,}}

\def\BCL{Barbour, Chen \& Loh}
\def\ABT{Arratia, Barbour \& Tavar\'e}
\def\BHJ{Barbour, Holst \& Janson}

\def\Blb{\left\{}
\def\Brb{\right\}}

\def\giv{\,|\,}
\def\Giv{\,\Big|\,}
\def\ep{\hfill$\Box$}
\def\bone{{\bf 1}}
\def\non{\nonumber}
\def\tI{{\widetilde I}}
\def\tS{{\widetilde S}}
\def\bh{{\bar h}}
\def\tp{{\tilde p}}
\def\tP{{\widetilde P}}

\def\adb#1{{\bf #1}}
\def\sodmi{\sum_{s=0}^{d-1}}
\def\sie{\sum_{i\in E}}
\def\sii{\sum_{i\ge1}}
\def\Le{\ \le\ }
\def\TP{{\rm TP\,}}
\def\Ord#1{O\Bl #1 \Br}
\newcommand{\eqs}{\begin{eqnarray*}}
\newcommand{\ens}{\end{eqnarray*}}

\newcommand{\eqa}{\begin{eqnarray}}
\newcommand{\ena}{\end{eqnarray}}
\newcommand{\eq}{\begin{equation}}
\newcommand{\en}{\end{equation}}

\def\numberlikeadb{\global\def\theequation{\thesection.\arabic{equation}}}
\numberlikeadb
\newtheorem{theorem}{Theorem}[section]
\newtheorem{lemma}[theorem]{Lemma}
\newtheorem{corollary}[theorem]{Corollary}
\newtheorem{proposition}[theorem]{Proposition}
\newtheorem{example}[theorem]{Example}

\def\pr{{\mathbb P}}
\def\ex{{\mathbb E}}
\def\re{{\mathbb R}}
\def\integ{{\mathbb Z}}
\def\ZZ{\integ}
\def\Ceka{\v Cekan\-avi\v cius}
\def\sjom{\sum_{j=0}^m}
\def\uo{^{(0)}}
\def\ul{^{(l)}}
\def\gg{{\mathcal G}}
\def\bh{{\bar h}}
\def\bX{{\widebar X}}
\def\bU{{\widebar U}}
\def\BX{Barbour \&~Xia}
\def\BC{Barbour \&  \Ceka}
\def\sro{\sum_{r\ge0}}
\def\srz{\sum_{r\in\integ}}
\def\ssz{\sum_{s\in\integ}}
\def\diff{\lfloor \m-\s^2 \rfloor}
\def\diffi{\lfloor \m_1-\s_1^2 \rfloor}
\def\difft{\lfloor \m_2-\s_2^2 \rfloor}
\def\sid{\s_1^2 + \d_1}
\def\std{\s_2^2 + \d_2}
\def\hX{{\widehat X}}
\def\hY{{\widehat Y}}
\def\hW{{\widehat W}}
\def\hQ{{\widehat Q}}
\def\hS{{\widehat S}}
\def\Ge{{\rm Ge\,}}
\def\sjn{\sum_{j=1}^n}
\def\bone{{\bf 1}}
\def\ixo{\bone_{\xx_0}}
\def\ixoc{{\bf 1}_{\xx_0^c}}
\def\aaomi{\aaa_0^{-1}}
\def\uu{{\cal U}}
\def\Eq{\ =\ }
\def\bb{{\cal B}}
\def\aaa{{\cal A}}
\def\sjo{\sum_{j\ge0}}
\def\dhh{d_{\hh}}
\def\dff{d_{\bff}}
\def\ngg{\|_{\gg}}
\def\bff{{\overline\ff}}
\def\nsup{\|_\infty}
\def\thalf{\tfrac12}
\def\derivbnd{\sqrt{2\pi}}
\def\sli{\sum_{l\ge1}}
\def\nn{{\cal N}}
\def\hh{{\cal H}}
\def\ui{^{(1)}}
\def\ut{^{(2)}}
\def\um{^{(m)}}
\def\ul{^{(l)}}
\def\ud{^{(\cdot)}}
\def\tih{{\tilde\h}}
\def\tL{{\widetilde L}}
\def\hL{{L'}}
\def\ddt{{t'}}
\def\hdd{{\h'}}
\def\tihd{{\tilde{\h'}}}

\def\ignore#1{}
\def\xxo{{\bf X}}
\def\LL{{\overline L}}
\def\tl{{\tilde\l}}
\def\Ge{\ \ge\ }
\def\Lip{{\rm Lip}}
\def\dpt{d_{{\rm pt}}}
\def\sjz{\sum_{j\in\integ}}
\def\prax{Proposition~\ref{appx}}
\def\sko{\sum_{k\ge0}}
\def\slo{\sum_{l\ge0}}
\def\sjk{\sum_{j=1}^k}
\def\hP{{\widehat P}}
\def\sjj{\sum_{j\in J}}
\def\xue{X^u}
\def\xle{X^{l,\e}}
\def\bx{{\bar x}}
\def\by{{\bar y}}
\def\RC{{\bf RC\,}}

\def\uii{^{(i)}}
\def\twi{{\widetilde W}\uii}
\def\uil{u_{il}}
\def\vil{v_{il}}
\def\sln{\sum_{l=1}^n}
\def\BP{{\rm BP}}
\def\ignore#1{}

\title{Coupling a branching process to an infinite dimensional epidemic process}
\author
{A.\ D.\ Barbour\thanks{Angewandte Mathematik,
Winterthurerstrasse~190,
CH--8057 Z\"URICH, Switzerland: {\tt a.d.barbour@math.uzh.ch}
}
\\
Universit\"at Z\"urich
}
\date{}
\maketitle

\vglue-1cm
\begin{abstract}
Branching process approximation to the initial stages of an epidemic
process has been used since the 1950's as a technique for providing 
stochastic counterparts to deterministic epidemic threshold theorems.
One way of describing the approximation is to construct both 
branching and epidemic processes on the same probability space, in
such a way that their paths coincide for as long as possible.  In
this paper, it is shown, in the context of a Markovian model of parasitic
infection, that coincidence can be achieved with asymptotically high 
probability until~$M_N$ infections have occurred, as long as
$M_N = o(N^{2/3})$, where~$N$ denotes the total number of hosts.
\end{abstract}

%\vskip10pt

\section{Introduction}\label{Intro}
 \setcounter{equation}{0}
%\subsection{The SIR model}

The classical law of large numbers and central limit theorem have process
analogues for many Markovian models arising in population ecology.
The law of large numbers is replaced by a deterministic process, obtained by
solving an appropriate system of ordinary or partial differential equations, and
the central limit theorem is replaced by a diffusion approximation around the
deterministic limit.  For many techniques and examples concerning such
density dependent Markov population processes, see Kurtz (1976, 1981).

In the context of invasion biology, when the central question is whether
the introduction of a small number of individuals of a species can lead to
its becoming established in a new habitat, these large population approximations
are no longer appropriate.  The more natural process approximations, at least 
if spatial restrictions on mixing are not critical in such small populations,
are now branching processes. These were introduced, in the context of 
epidemic theory, by Whittle~(1955), Kendall~(1956) and Bartlett~(1956, p.~129);
here, infected individuals play the part of the invading species, and those
that are infected by an individual correspond to an individual's `offspring'.
  
When considering the development of a single species as a branching process, 
the biological quantity~$R_0$, the lifetime mean number of offspring of a 
single individual
when unhampered by competition from others of the same species,
is just the mean offspring number of the corresponding Galton--Watson process.
The branching process criticality theorem then corresponds to the 
biological meta-theorem,
that an invading population can only become established if its $R_0$ (in
the context that it experiences upon invasion) exceeds~$1$.  For models
involving more species, the analogy is to multitype branching processes,
and the dominant eigenvalue of the mean matrix of the branching
process has a corresponding interpretation in the biological context.
For more detailed discussions of such issues, see Heesterbeek~(1992)
and Diekmann \& Heesterbeek~(2000, Section~5.7).

Whittle~(1955) was able to justify his birth and death approximation to the
early stages of the Markovian SIR-epidemic, and hence his formula for the
probability of a large epidemic occurring, by sandwiching the epidemic
process, during its initial stages,  between two birth and death processes 
with slightly differing transition rates.  This can be interpreted in terms
of a pathwise comparison of processes. Ball~(1983) and
Ball and Donnelly~(1995) went rather further, using a coupling
argument to link the epidemic process with an approximating branching process
on one and the same probability space, in such a way that the paths of the
two processes are identical for a certain (random) length of time.  
In particular, they showed that the 
total variation distance between the distributions of the paths of the 
branching and epidemic
processes is small, up to the time at which $M=M_N$ infections have taken place,
for any choice $M_N = o(\sqrt N)$.  They also suggest that this range
of~$M_N$ cannot be extended.

The coupling used by Ball and Donnelly is simple and natural, and it is
somewhat surprising that accurate coupling is in fact possible, for
some epidemic processes, over rather longer time intervals than they
had supposed possible.  This was first established by Barbour and 
Utev~(2004), in the context of the discrete time Reed--Frost epidemic process.
They showed that the branching process approximation to the path
distribution actually has asymptotically small error in total
variation for all choices of~$M_N = o(N^{2/3})$.  
The essence of their
argument lay in examining the likelihood ratio of the two processes along
paths of given length, and showing that it was typically close to~$1$.  In
this paper, we show that similar arguments can also be applied to some continuous
time models. We take as example the infinite dimensional BK-model, 
introduced in Barbour and Kafetzaki~(1993) and subsequently generalized
by Luchsinger~(2002a,b), for describing the transmission
of the parasitic disease schistosomiasis.
\ignore{
A first theorem shows that the path distributions can
be matched asymptotically accurately
over paths of length~$m_N$ transitions, as long as $m_N = o(N^{2/3})$.
The second, stronger theorem proves a similar result for paths in which
at most~$M_N$ infections have taken place.
  In this model, for some
ranges of the parameters, $M_N$ may be significantly smaller than~$m_N$.
}

\section{The BK-model}\label{model}
 \setcounter{equation}{0}
In the BK-model, $N$ hosts are infected by parasites, with $X_j^N(t)$
hosts having~$j$ parasites at time~$t$, for $j\in\ZZ_+$ and $t\ge0$.
The process evolves as a Markov jump process~$X^N$ in continuous time on the
set $\xx := \{(\xi_j \in \ZZ_+,\,j\ge0)\colon\,\sjo\xi_j = N\}$, with
transition rates given by
\eqs
  \xi &\to& \xi + e(j-1) - e(j) \quad\mbox{at rate}\quad j\m \xi_j,\quad j\ge1;\\
  \xi &\to& \xi + e(j) - e(0) \quad\mbox{at rate}\quad \l\xi_0\sli (\xi_l/N)p_{lj},
   \quad j\ge1,
\ens
for any $\xi\in\xx$,
where $e(j)$ denotes the unit vector in the $j$-th coordinate.  The first of the
transitions models the death of a parasite in one of the~$\xi_j$ hosts currently
carrying~$j$ parasites, the parasites being assumed to have independent exponentially
distributed lifetimes with mean~$1/\m$.  The second transition models infection.
Only currently uninfected hosts can be newly infected, and each makes contacts
that could potentially lead to
infection at rate~$\l$, the chance of such a contact being made with
a host carrying~$l$ parasites being $(\xi_l/N)$ (homogeneous mixing of hosts).
If there is such a contact between an uninfected host and an~$l$-host, then~$j$
parasites are established in the previously uninfected host with probability~$p_{lj}$;
in the BK-model, it is supposed that $p_{lj} = \pr[U_l = j]$, for $U_l := \sum_{i=1}^l
Y_i$, where the~$Y_i$ are independent and identically distributed random variables
with mean~$\th$ and finite variance, implying that each of the~$l$ parasites 
transmits on average~$\th$ infective stages to the newly infected host at an infectious
contact, independently of the others.

For a disease such as schistosomiasis, infection is actually indirect, and involves
a host infecting suitable aquatic snails and these snails subsequently passing 
infection to other hosts. Thus the BK-model does not seem at first sight to be
at all realistic.  However, it can be thought of as an extreme case of a model
incorporating features of the transmission process that were not present
in many of the previous models:  infection by parasites in groups, rather than
singly, immunity in the definitive host (here, in the form of perfect concomitant 
immunity), explicit incorporation of the parasite burdens of individual hosts.

The model that results is interesting for a number of reasons.  The first is that,
although it is rather complicated, it is still simple enough for some analytic
conclusions to be reached.  For instance, it can be shown that the model has a
`law of large numbers' approximation for large~$N$, in the form of the solution
to an infinite system of differential equations, whose components approximate the
proportions of hosts with different numbers of parasites. If~$\th > e$, this
differential equation system has no (endemic) equilibrium solution that yields
a finite mean number of parasites per host.  In practice, the distribution of
parasites among hosts is observed to be extremely irregular, so that such
behaviour is very encouraging: most earlier models have tacitly predicted Poisson--like
distributions, which are far from realistic, and those that have tried to account
for the over-dispersed distributions observed have imposed a specific form 
for the distribution without
proposing any mechanism that might generate it.  Another feature is that, if $\th<e$, 
there is exactly one equilibrium distribution of the differential equation system
that has finite mean number of parasites per host, and that, in this
equilibrium, the distribution of the
number of parasites per host, conditional on the host being infected, depends
only on the value of~$\th$, and not on $\l$ or~$\m$.

For the purposes of this paper, it is the behaviour when few hosts are infected
that is of primary relevance, with interest centering on questions such as the
probability that the introduction of a single infected host can cause 
endemic infection to become established.  These are the kinds of problem
that can be addressed by way of a branching process approximation.
Here, we begin by proving an error bound for the approximation (Theorem~\ref{BPA-T1})
that is asymptotically 
valid in total variation for paths of length~$o(N^{2/3})$ transitions 
as $N\to \infty$. The branching process in turn yields a criticality theorem, 
which, to a close approximation,
describes whether or not endemic equilibrium is possible in the BK-model. 
 
However, the approximating
branching process --- a Markov branching process with countably infinitely many
types --- itself displays unexpected critical behavour.  If $\th\le e$, the
branching process is super-critical, in the sense of having positive probability 
of growing indefinitely, if and only if $\l\th/\m > 1$.
The quantity $\l\th/\m$ has an immediate interpretation,
being the lifetime average number of offspring of a single parasite, where offspring
is interpreted in terms of parasites successfully passed on to other hosts, 
and is therefore precisely the biological quantity~$R_0$, as seen from
the parasites' viewpoint.  Its appearance as the criticality parameter is 
therefore exactly what one would expect.  However,
if $\th > e$,  the criticality parameter is
$\l e \log\th/\m$, a fact that is much more difficult to interpret.

Another feature of the model is that the mean number of parasites
develops in time with exponential rate~$\l\th-\m$, whereas, if $\th>e$, a 
super-critical process has a smaller exponential growth rate for the number 
of infected hosts. Thus, in such circumstances, the mean number of parasites 
per host increases ever
faster.  As a result, because deaths of parasites are counted
as transitions, paths containing~$m_N$ transitions may contain many fewer
infections --- roughly speaking, one may well have only
$M_N \approx m_N^\a$ infections, for some~$\a < 1$.
For such choices of the parameters, this makes the above theorem unsuitable for 
direct comparison with the results of Ball and Donnelly~(1995).  We
therefore prove a second error bound in Theorem~\ref{BPA-T2}, which is  
expressed in terms of the asymptotics
of~$M_N$.  Its proof turns out to be a relatively simple adaptation
of that of Theorem~\ref{BPA-T1}. We conclude with Theorem~\ref{local},
which establishes a rather stronger local statement, showing that the ratios
of the likelihoods under the two models of paths containing at 
most~$M$ infections typically differ from~$1$ by more than order
$O\Bl(M^{2/3}/N)\sqrt{\log(N/M^{2/3})}\Br$ with asymptotically
negligible probability.

\section{Total variation approximation}\label{BPA}
 \setcounter{equation}{0}
The Markov branching process $X := (X_j(\cdot),\,j\ge1)$ that approximates the BK-model
is obtained from the process  $X^N$ by ignoring the $0$-component, taking the 
countable set 
$\xx^* := \{(\xi_j \in \ZZ_+,\,j\ge1)\colon\,\sji\xi_j < \infty\}$ as
state space, and modifying the transition rates to
\eqs
  \xi &\to& \xi + e(j-1) - e(j) \quad\mbox{at rate}\quad j\m \xi_j,\quad j\ge2;\\
  \xi &\to& \xi  - e(1) \quad\mbox{at rate}\quad \m \xi_1,\\ 
  \xi &\to& \xi + e(j) \quad\mbox{at rate}\quad \l\sli \xi_l p_{lj},
   \quad j\ge1,
\ens   
for $\xi \in \xx^*$.
These rates are identical with those for $X^N$, except that, in the infection
transition, the factor $\xi_0/N = 1 - S(\xi)/N$ is replaced by~$1$, 
where $S(\xi) := \sji\xi_j$.  This represents the fact that, in the branching
approximation, the total proportion of infected hosts is considered to be vanishingly
small.  Clearly, this should make little difference to individual transitions
if $S(\xi) \ll N$.  The main result of this paper is to show that it makes little
difference even for the distribution of whole path segments, considered as
paths in~$\xx^*$, provided that the
number of transitions~$m$ in the segment and the initial 
state~$\xi\uo \in \xx^*$ are such
that $m + S(\xi\uo) \ll N^{2/3}$.  We denote such a path by
$\{(\xi\ul,t\ul),\,0\le l\le m\}$, where $t\uo := 0$, and we let~$\ff_m$ denote 
the Borel $\s$-algebra of events generated by these paths.  To avoid trivial
exceptions caused by paths that are absorbed in~$0 \in \xx^*$ never making
further jumps, we suppose that both processes, when in state $0$, make `jumps'
to state~$0$ at unit rate.

\begin{theorem}\label{BPA-T1}
Suppose that $\xi\uo\in \xx^*$, $N\ge2$ and~$m$ are such that
$S_m \sqrt m < N$, where $S_m := m + S(\xi\uo)$.
Then, for any $A \in \ff_m$, we have
\[
   |\pr[X^{N*} \in A] - \pr[X \in A]| 
       \Le 8\frac{S_m\sqrt m}N,
\]
% provided that $m \ge (1350/13)\log N$,
where $X^{N*}$ denotes the process~$X^N$ without the zero coordinate.
% For smaller~$m$, the inequality is true if the factor~$20$ is replaced by~$85$.
\end{theorem}

\proof
For $\xi\in\xx^*$ with $1 \le S(\xi) \le N$, write
\eqs
    \L^N(\xi) &:=& \l(1-S(\xi)/N)\sli \xi_l (1 - p_{l0}),\qquad 
    \rho^N(\xi) \Eq \m\sji j\xi_j + \L^N(\xi);\\
		\L(\xi) &:=& \l \sli \xi_l (1 - p_{l0}),\qquad 
    \rho(\xi) \Eq  \m\sji j\xi_j + \L(\xi).
\ens
The quantities $\rho^N(\xi)$ and~$\rho(\xi)$ respectively 
denote the overall jump rates of the processes $X^N$ and~$X$ in state~$\xi$,
$\L^N(\xi)$ and~$\L(\xi)$ the overall infection rates; for
$\xi = 0 \in \xx^*$, we set $\rho^N(0) = \rho(0) = \L^N(0) = \L(0) = 1$.
Suppose that $m+S(\xi\uo) \le N$.
Then, for a path with~$m$ transitions starting in~$\xi\uo$ at time~$0$
and then passing through the sequence of states $(\xi\ul,\,1\le l\le m)$
at times $t\ui < t\ut < \cdots < t\um$, the likelihood ratio $d\pr_{X^N}/d\pr_X$
evaluated at such a path is just
\eqa
  \lefteqn{L_m^N \ :=\ L_m^N(\xi\ud,t\ud)} \non\\
   &:=& \prod_{l=0}^{m-1} \Blb \exp\{-(\L^N(\xi\ul) - \L(\xi\ul))(t^{(l+1)}-t\ul)\} 
%     \frac{\L(\xi\ul)}{\L^N(\xi\ul)}
		 \,\Bl 1 - \frac{S(\xi\ul)}{N} \Br^{u_l} \Brb, \label{LR-I}
\ena
where
\[
    u_l \Eq 
\begin{cases}
               1 &\mbox{if}\quad S(\xi^{(l+1)}) = S(\xi\ul) + 1;\\
               0 &\mbox{otherwise}.
\end{cases}							 
\]
Hence we have
\[
   L_{l+1}^N \Eq L_l^N(1+\h_{l1}^N)(1-\h_{l2}^N),  %(1+\h_{l3}^N),
\]
with
\eqs
   \h_{l1}^N &:=& \exp\{-(\L^N(\xi\ul) - \L(\xi\ul))(t^{(l+1)}-t\ul)\} - 1,\\
%   \h_{l2}^N &:=& \frac{\L(\xi\ul)}{\L^N(\xi\ul)} - 1
\ens
and
\[
    \h_{l2}^N \ :=\  \frac{S(\xi\ul)}{N} \bone\{u_l=1\}.
\]
Note that each of these quantities is zero if $\xi\ul = 0 \in \xx^*$.

The inequality
\eq\label{eta-3}
  0 \Le \h_{l2}^N \Le S(\xi\ul)/N
\en
is immediate.
Then, from the definitions of $\L^N$ and~$\L$, it follows directly that
\[
   \L(\xi) - \L^N(\xi) \Eq \l \frac{S(\xi)}N \sji \xi_j(1-p_{j0}),
\]
implying that
\eq\label{L-ratio-1}
   |\L^N(\xi)/\L(\xi) - 1| \Le S(\xi)/N.
\en
\ignore{
Hence, provided also that~$\xi\ul$ is such that
\eq\label{star}
    S(\xi\ul)/N \Le \frac12,
\en
it follows that
\eq\label{eta-2}
  |\h_{l2}^N| \Le 2S(\xi\ul)/N.
\en
} 
Furthermore, if
\eq\label{2-star}
  \tih_{l1}^N \ :=\ |\L^N(\xi\ul)/\L(\xi\ul) - 1| \L(\xi\ul)(t^{(l+1)}-t\ul) \Le 1,
\en
we also have
\eq\label{eta-1}
  |\h_{l1}^N| \Le 2|\tih_{l1}^N| \Le 2\{S(\xi\ul)/N\} e_l,
\en
where
\[
   e_l \ :=\ \rho(\xi\ul)(t^{(l+1)}-t\ul).
\]
Hence, if %\Ref{star} and
\Ref{2-star} is satisfied, it follows that
\eq\label{cross}
  |L_{l+1}^N - L_l^N| \Le N^{-1}S(\xi\ul)\{1 + 2e_l\} L_l^N.
\en

Now suppose that $(X\ul,\,l\ge0)$ is a path resulting
from a realization of the process~$X$ starting with $X(0) = \xi\uo$, 
and that $(T\ul,\,l\ge1)$ are the
corresponding jump times: set~$T\uo = 0$.  Then, defining 
\eq\label{E-defn}
   E_{l+1} \ :=\ \rho(X\ul)(T^{(l+1)}-T\ul),\qquad l\ge0,
\en
we note that $\law(E_{l+1}\giv \ff_l)$ is the standard exponential
distribution for each~$l$.  Furthermore, the process 
$\{L_l^N(X\ud,T\ud),\, l\ge0\}$
is a non-negative martingale with $L_0^N = 1$ a.s., and so is the
stopped version $\tL_l^N := L_{l\wedge\t_1^N\wedge\t_2^N}^N(X\ud,T\ud)$, 
where
\eqs
   \t_1^N &:=& \inf\{l\ge0\colon\, E_{l+1} > N/S(\xi\ul)\},\\
   \t_2^N &:=& \inf\{l\ge1\colon\, L_l^N > 2\}.
\ens
Note also that $S(X\ul) \le S(\xi\uo) + l \le S_m$
% , so that~\Ref{star} is satisfied 
for all $0\le l\le m$. % because $S_m  < N/2$.
%Since~\Ref{star} is satisfied for all $0\le l\le m$, i
Now it follows 
from \Ref{L-ratio-1} and~\Ref{E-defn} that~\Ref{2-star} is satisfied for all
$0 \le l < \t_1^N$.  
Hence, from~\Ref{cross} and the definition of~$\t_2^N$, it follows that
\[
  |\tL_{l+1}^N - \tL_l^N| 
	     \Le 2S_mN^{-1}(1 + 2E_l)\quad \mbox{for all}\ 0 \le l < m.
\]		
% We are thus in a position to apply Lemma~\ref{app} to the martingale~$\tL_l$,
% with $n=m$ and $2a=b=4S_m/N$.  %, and any~$y$ such that
Thus, because~$\tL_l^N$ is a martingale, it follows that
\eq\label{variance}
  \ex(\tL_m^N-1)^2 \Le \frac{4mS_m^2}{N^2}\,\ex(1 + 2E_1)^2 
	           \Eq \frac{52mS_m^2}{N^2}.
\en

Now, for any $A \in \ff_m$, we have
\eqa
	 \lefteqn{\pr[X \in A] - \pr[X^{N*} \in A] \Eq \ex\{(1- L_m^N) \bone\{X\in A\} \} 
	 \Le \ex\{(1- L_m^N)^+\}} \non\\
	 &\le& \pr[\t_1^N \le m] + \pr[\{\t_2^N \le m\}\cap\{\t_1^N > m\}]
%	     + \int_0^1 \pr[1 - \tL_m^N > y]\,dy.\phantom{HH} 
       + \ex(1 - \tL_m^N)^+.
			 \label{difference}
\ena
{}From the definition of~$\t_1^N$, it is immediate that
\eq\label{t1-bnd}
   \pr[\t_1^N \le m] \Le me^{-N/S_m} \Le 4e^{-2}\,\frac{S_m\sqrt m}N
\en
if $S_m\sqrt m \le N$.	
Then we have
\eq\label{t2-bnd}
    \pr[\{\t_2^N \le m\}\cap\{\t_1^N > m\}] \Le \pr[\tL_m^N - 1 > 1] 
		       \Le \ex(\tL_m^N - 1)^+.
\en
Finally, we have
\eqs
   \ex(1 - \tL_m^N)^+ + \ex(\tL_m^N - 1)^+ &=& \ex|1 - \tL_m^N| \\
	 &\le& \sqrt{\ex(1 - \tL_m^N)^2} \Le 2\sqrt{13}S_m\sqrt m / N.
\ens
It remains to note that $4e^{-2} + 2\sqrt{13} < 8$.
\ep

\medskip		 
In general, the bound given in the theorem provides useful information
as long as $S_m\sqrt m \ll N$.  In asymptotic terms, 
for fixed~$\xi\uo$, this allows paths of lengths $m_N = o(N^{2/3})$ as 
$N\to\infty$, with an error bound of order $O(N^{3\g/2-1})$ 
if $m_N \sim N^\g$ for some $\g < 2/3$.  

If the Ball and Donnelly~(1995) coupling is used to obtain error bounds,
the resulting order~$O(N^{2\g-1})$, if $\xi\uo$ is fixed
and $M_N \sim N^\g$, is at first sight not as sharp.  However, there is an 
important difference between the two results:  the theorem
above has~$m_N$, the total number of transitions, in the error bound,
whereas the Ball and Donnelly coupling leads to an error expressed in terms 
of~$M_N$, the number of births or infections.
Now the total number of transitions includes all the parasite deaths, and 
if the mean number of parasites per host grows fast, as may be the case when 
$\th > e$, $m_N$ may be substantially bigger than~$M_N$.  Thus Theorem~\ref{BPA-T1}
is not strong enough to yield an obvious improvement.  For this reason,
we now bound the discrepancies in the likelihood ratio more carefully,
basing the argument explicitly on the sequence of infection events.
To this end, we let~$\hh_l$ denote the Borel $\s$-algebra of events generated 
by paths containing exactly~$l$ infection events; as before, to avoid trivial
exceptions caused by paths that are absorbed in~$0 \in \xx^*$ having no
further infections, we suppose that both processes, when in state $0$, 
create `pseudoinfections' at unit rate.

\begin{theorem}\label{BPA-T2}
Suppose that $\xi\uo\in \xx^*$, $N\ge2$ and~$M$ are such that
$S_M \sqrt M \le N$, where $S_M := M + S(\xi\uo)$.
Then, for any $A \in \hh_M$, we have
\[
   |\pr[X^{N*} \in A] - \pr[X \in A]| 
       \Le 8\frac{S_M\sqrt M}N.
\]
% provided that $M \ge (1350/13)\log N$.
% where $X^{N*}$ denotes the process~$X^N$ without the zero coordinate.
% For smaller~$M$, the inequality is true if the factor~$20$ is replaced by~$85$.
\end{theorem}

\proof
The likelihood ratio at a path~$\xi(\cdot)$ in~$\hh_M$ can be written, 
using~\Ref{LR-I}, in the form
\eqa
  \lefteqn{\hL_M^N \ :=\ \hL_M^N(\xi(u),\,0\le u\le \s_M)} \non\\
   &:=& \prod_{l=1}^{M} \Blb \exp\Bl 
	    -\int_{\s_{l-1}}^{\s_l} (\L^N(\xi(u)) - \L(\xi(u)))\,du \Br 
%     \frac{\L(\xi\ul)}{\L^N(\xi\ul)}
		 \,\Bl 1 - \frac{S(\xi(\s_l))}{N} \Br \Brb,\phantom{HH} \label{LR-II}
\ena
where $0 = \s_0 < \s_1 < \cdots$ denote the times of infection transitions.
Hence, very much as before, we have
\[
   \hL_{l}^N \Eq \hL_{l-1}^N(1+\hdd_{l1}^N)(1-\hdd_{l2}^N),  %(1+\h_{l3}^N),
\]
with
\eqs
   \hdd_{l1}^N &:=& \exp\Bl 
	      -\int_{\s_{l-1}}^{\s_l} (\L^N(\xi(u)) - \L(\xi(u)))\,du \Br - 1,\\
%   \h_{l2}^N &:=& \frac{\L(\xi\ul)}{\L^N(\xi\ul)} - 1
\ens
and
\[
    \hdd_{l2}^N \ :=\  \frac{S(\xi(\s_l))}{N},
\]
each of these quantities being zero if $\xi\ul = 0 \in \xx^*$.

{}From~\Ref{L-ratio-1}, setting
\eq\label{tih-def-2}
  \tihd_{l1}^N \ :=\  
	    \int_{\s_{l-1}}^{\s_l} \{1 - \L^N(\xi(u))/\L(\xi(u))\}\L(\xi(u))\,du
\en
we have
\eq\label{eta-1-2}
  0 \Le \hdd_{l1}^N \Le 2\tihd_{l1}^N \Le 2\{S(\xi(\s_{l-1}))/N\} e_l,
\en
whenever $\tihd_{l1}^N \le 1$, where 
\[
   e_l \ :=\ \int_{\s_{l-1}}^{\s_l} \L(\xi(u))\,du.
\]
Hence, noting that $S(\xi(u)) \le S_M$ for $0 \le u \le \s_M$, if %\Ref{star} and
\eq\label{2-star-2}
   e_l \Eq \int_{\s_{l-1}}^{\s_l} \L(\xi(u))\,du \Le N/S_M
\en	 
is satisfied, it follows that
\eq\label{cross-2}
  |\hL_{l}^N - \hL_{l-1}^N| \Le N^{-1}S_M\{1 + 2e_l\} \hL_{l-1}^N, \quad 1\le l\le M.
\en

Now suppose that $(X(u),\,u\ge0)$ is a path resulting
from a realization of the process~$X$ starting with $X(0) = \xi\uo$, 
and that $(\s_l,\,l\ge1)$ are the
corresponding times of births (infection transitions): set~$\s_0 = 0$.  
Then, defining 
\eq\label{E-defn-2}
   E'_{l} \ :=\ \int_{\s_{l-1}}^{\s_l} \L(X(u))\,du,\qquad l\ge1,
\en
we note that $\law(E'_{l+1}\giv \hh_l)$ is the standard exponential
distribution for each~$l$.  We now argue as before using the 
likelihood ratio martingales $\{\hL_l^N(X(\cdot)),\, l\ge0\}$
and $\tilde\hL_l^N := \hL_{l\wedge\ddt_1^N\wedge\ddt_2^N}^N(X(\cdot))$,  
where
\eqs
   \ddt_1^N &:=& \inf\{l\ge0\colon\, E'_{l+1} > N/S_M\},\\
   \ddt_2^N &:=& \inf\{l\ge1\colon\, \hL_l^N > 2\}.
\ens
Since~\Ref{2-star-2} is satisfied for all $0 \le l < \ddt_1^N$, it follows 
from~\Ref{cross-2} and the definition of~$\ddt_2^N$, that
\eq\label{hL-diff}
  |\tilde\hL_{l+1}^N - \tilde\hL_l^N| 
	     \Le 2S_M N^{-1}(1 + 2E'_{l+1})\quad \mbox{for all}\ 0\le l < M.
\en		
The remaining argument is now exactly as before, using % Lemma~\ref{app} and
the martingale~$\tilde\hL_l$ to compare the probabilities 
$\pr[X^{N*} \in A]$ and~$\pr[X \in A]$ for $A \in \hh_M$.
\ep 

\medskip
Thus Theorem~\ref{BPA-T2} yields bounds of order $O(N^{3\g/2-1})$, improving on the 
rate obtained using the Ball and Donnelly~(1995) coupling, if $\xi\uo$ is fixed
and $M_N \sim N^\g$ for $\g < 3/2$, where $M_N$ denotes the number of 
infection transitions.

The new argument exploits the fact that the life histories of individuals
infected with a given number of parasites have identical distributions in
both models, except for the infection events, so that the likelihood ratio
is correspondingly simpler.  The key element is then that the {\em difference} in
infection rates between the two models is sufficiently small compared to the
infection rate itself.  The argument in Theorem~\ref{BPA-T1} is less precise
largely because, if the number of parasites is large, the bound~\Ref{eta-1}
is rather pessimistic, since a potentially small factor $\L(\xi)/\rho(\xi)$
is not being exploited.

\section{Local approximation}\label{BPA2}
 \setcounter{equation}{0}
It was argued in Barbour and Utev~(2004) that total variation approximation
is not necessarily the best measure of closeness, if statistical applications
involving likelihoods are to be justified.  It is much more natural to want
to have local approximations, which ensure that the ratio of actual and
approximate likelihood is very close to~$1$, except possibly on a set of very
small probability.  As a result, they defined a measure of relative closeness:
probability measures $P$ and~$Q$ on~$\ff$ are said to be $\e$-relatively close with
tolerance~$\h$, $\RC(\e,\h)$ for short, if there exists a set~$R\in\ff$ such
that
\[
   P(R^c) \Le \h,\quad Q(R^c) \Le \h,\quad \sup_{x\in R}|\log((dP/dQ)(x))| \Le \e.
\]
In this section, we show that the branching process approximation of the
previous sections is indeed relatively close, as long as $M_N \ll N^{2/3}$.
	 
We begin with a minor modification of the bounded differences
lemma for martingales.

\begin{lemma}\label{app}
If $(L_n,\gg_n,\,n\ge0)$ is a martingale, and if
\[
   |L_{n+1} - L_n| \Le a + bE_{n+1} \quad \mbox{for each}\ n\ge0,
\]
where $\law(E_{n+1}\giv \gg_n)$ is the standard exponential distribution
$\exp(1)$ for each $n\ge1$, then
\[
  (1):\quad\ 
	  \max\{\pr[L_n - L_0 \le -y], \pr[L_n - L_0 \ge y]\} 
		      \Le \exp\Blb \frac{-3y^2}{8n\{(a+b)^2 + b^2\}} \Brb
\]
for all 
\[
   0 \Le y \Le \frac{4n}3 \e_0\{(a+b)^2 + b^2\}/\max(a,b),
\]
where $\e_0 > 1/15$ is the constant defined by $e^{\e_0}(1-\e_0)^{-3} = 4/3$.
Furthermore, for all $y\ge0$,
\[
   (2):\quad 
	   \max\{\pr[L_n - L_0 \le -y], \pr[L_n - L_0 \ge y]\} \Le
		     \exp\Blb \frac{-y}{15\max(a,b)\sqrt n} + \frac2{135} \Brb.
\]				 
\end{lemma}

\proof
If~$X$ is any random variable with $\ex X = 0$ and $|X| \le a + bE$, where
$E \sim \exp(1)$, then it follows that, for any $\th > 0$,
\eqs
  \ex\{e^{\th X}\} &\le& \ex\Blb 1 + \th X + \half \th^2 X^2 e^{\th|X|}\Brb\\
  &\le& 1 + \half \th^2 \ex\Blb (a+bE)^2 e^{\th(a+bE)} \Brb \\
  &=& 1 + \half \th^2 e^{a\th}\Blb \frac{a^2}{1-b\th} + \frac{2ab}{(1-b\th)^2}
      + \frac{2b^2}{(1-b\th)^3} \Brb \\
  &\le&  \exp\Blb \frac23 \th^2 \{(a+b)^2 + b^2\} \Brb,
\ens
as long as $\th\max(a,b) \le \e_0$, with $\e_0$ defined as above.
Hence, for any $n\ge1$,
\eqs
  \ex\Blb e^{\th(L_n - L_0)} \giv \gg_{n-1}\Brb &=& 
    e^{\th(L_{n-1}-L_0)} \ex\Blb e^{\th(L_n - L_{n-1})} \giv \gg_{n-1}\Brb \\
  &\le& \exp\Blb \frac23 \th^2\{(a+b)^2 + b^2\} \Brb e^{\th(L_{n-1}-L_0)} ,
\ens
implying that
\[
  \ex\Blb e^{\th(L_n - L_0)} \Brb \Le \exp\Blb \frac23 \th^2\{(a+b)^2 + b^2\} \Brb 
      \ex\Blb e^{\th(L_{n-1} - L_0)} \Brb
\]
for all $n\ge 1$, and hence that
\[
  \ex\Blb e^{\th(L_n - L_0)} \Brb \Le \exp\Blb \frac23 n\th^2\{(a+b)^2 + b^2\} \Brb.
\]
Hence, for any $y\ge0$ and any $\th$ such that $\th\max(a,b) \le \e_0$, we have
\[
   \pr[L_n - L_0 \ge y] \Le \exp\Blb -y\th + \frac23 n\th^2\{(a+b)^2 + b^2\} \Brb.
\]
Now, if $y\max(a,b) \le (4\e_0/3)n\{(a+b)^2 + b^2\}$, we can take
\[
   \th \Eq \frac y{(4n/3)\{(a+b)^2 + b^2\}},
\]
to give
\[
    \pr[L_n - L_0 \ge y] \Le \exp\Blb \frac{-3y^2}{8n\{(a+b)^2 + b^2\}} \Brb.
\]
On the other hand, for all $y\ge0$ and $n\ge1$, we can choose 
$\th = 1/\{15\max(a,b)\sqrt n\}$, giving
\[
   \pr[L_n - L_0 \ge y] \Le \exp\Blb \frac{-y}{15\max(a,b)\sqrt n} + \frac2{135} \Brb.
\]
	 
The same arguments also cover $\pr[L_n - L_0 \le -y]$ for the corresponding
choices of~$y$, since the conditions of the theorem apply equally well to 
the martingale~$-L_n$.
\ep

\medskip
This lemma enables us to prove the following estimate of relative closeness.

\begin{theorem}\label{local}
Suppose that $\xi\uo\in \xx^*$, $N\ge2$ and~$M$ are such that
$\ps(M,N) := S_M\sqrt M/N \le 1$, where $S_M := M + S(\xi\uo)$.  Then, with
respect to paths in~$\hh_M$, the processes~$X^{N*}$ and~$X$ are
$\RC(\e_{M,N}^r,\h_{M,N}^r)$ relatively close for any choice of $r \ge 1$, where
\eqs
   \e_{M,N}^r &:=& C_r \ps(M,N) \sqrt{\log(1/\ps(M,N))};\\
	 \h_{M,N}^r &:=&  2\ps(M,N)^r + e^{2/135}\exp\{-1/60\ps(M,N)\} + Me^{-N/S_M},
\ens	 
and $C_r := \sqrt{416r/3}$, provided that $M \ge (1/5)C_r^2 \log N$ and that
$\e_{M,N}^r \le 1$.	
\end{theorem}

\proof
It was shown in the proof of Theorem~\ref{BPA-T2} that the likelihoods of the
processes~$X^{N*}$ and~$X$ are close; here, we tighten the argument.
We start from~\Ref{hL-diff}, which states that 
\[
  |\tilde\hL_{l+1}^N - \tilde\hL_l^N| 
	     \Le 2S_M N^{-1}(1 + 2E'_{l+1})\quad \mbox{for all}\ 0\le l < M,
\]
where $\law(E'_{l+1}\giv \hh_l)$ is the standard exponential
distribution for each~$l$,
and from the observation that, by the definition of~$\tilde\hL_M^N$, we have 
$\tilde\hL_M^N = \hL_M^N$ as long as $\min\{\ddt_1^N,\ddt_2^N\} > M$.
Now it is immediate, as for~\Ref{t1-bnd}, that
\[
   \pr[\ddt_1^N \le M] \Le Me^{-N/S_M}.
\]	 
Then, from the definition of~$\ddt_2^N$, it follows that
\[
   \pr[\{\ddt_1^N > M\}\cap\{\ddt_2^N \le M\}] \Le \pr[\tilde\hL_M^N - 1 > 1].
\]
Hence, to establish the desired relative closeness, we take
\[
   R^c \ :=\ \{\min(\ddt_1^N,\ddt_2^N) \le M\} \cup \{|\tilde\hL_M- 1| > \e_{M,N}^r/2\},
\]
(here using the assumption that $\e_{M,N}^r \le 1$) 
and bound the probabilities $\pr[\tilde\hL_M^N - 1 > 1]$ 
and~$\pr[|\tilde\hL_M- 1| > \e_{M,N}^r/2]$ using Lemma~\ref{app} with $n = M$
and $2a = b = 4S_M/N$.

First, we use Lemma~\ref{app}\,(2) to give
\[
   \pr[\tilde\hL_M^N - 1 > 1] \Le \exp\{-N/(60S_M\sqrt M)\}e^{2/135}.
\]
Then we use Lemma~\ref{app}\,(1) to show that
\[
   \pr[|1 - \tilde\hL_M^N| > y] \Le 2\exp\{-3N^2y^2/(416M S_M^2)\},
\]
provided that
\[
   0 \Le y \Le \frac{4M}3\, \e_0 \,\frac{13S_M}N .
\]	 
Hence we can take $y=\e_{M,N}^r/2$ if 
\[
    C_r S_M\sqrt M \sqrt{\log N}/N \le 104MS_M/45N,
\]
and thus if $M \ge (1/5)C_r^2 \log N$, giving
\[
\hspace{0.2in}  \pr[|1 - \tilde\hL_M^N| > \e_{M,N}^r/2] 
	     \Le 2\Blb \frac{S_M\sqrt M}N \Brb^{(3/416)C_r^2} 
			   \Eq 2\Blb \frac{S_M\sqrt M}N \Brb^r. \hspace{0.2in}\Box
\]

\medskip
Thus asymptotic relative closeness of order 
$O\{ \ps(M,N) \sqrt{\log(1/\ps(M,N))}\}$
can be established with tolerance of arbitrarily small polynomial
order in $\ps(M,N) = S_M\sqrt M/N$.

\end{document}